# Enhanced Energy Management System with Corrective Transmission Switching Strategy— Part II: Results and Discussion

Xingpeng Li, *Student Member, IEEE* and Kory W. Hedman, *Member, IEEE*

**Abstract—** **This paper presents a novel procedure for energy management system (EMS) that can utilize the flexibility in transmission network in a practical way. With the proposed enhanced EMS procedure, the reliability benefits that are provided by corrective transmission switching (CTS) in real-time contingency analysis (RTCA) can be translated into significant cost savings in real-time security-constrained economic dispatch (RT SCED). Simulation results show the congestion cost with consideration of CTS is largely reduced as CTS can relieve potential post-contingency network violations. The effects of integrating CTS in existing EMS procedure on markets are also analyzed. In conclusion, this two-part paper shows that CTS can achieve substantial reliability benefits, as well as significant cost savings.**

## I. Introduction

P REVIOUS research in the literature has demonstrated that utilization of the flexibility in transmission networks offers a variety of benefits for real-time operations of electric power systems. The network topology is traditionally considered to be fixed. However, operators can reconfigure the network to achieve a specific target. Corrective transmission switching (CTS) can be considered as a practical strategy to utilize the flexibility in transmission networks.

CTS switches a transmission element out of service shortly after an outage to handle network violations. It can significantly reduce post-contingency network violations, which provides operators another option to handle contingencies. The CTS algorithms proposed in [1] are fast enough for real-time application, which is practical for actual implementation.

The main applications of energy management system (EMS) include real-time contingency analysis (RTCA) and real-time security-constrained economic dispatch (RT SCED). Although the authors' prior work [1]-[4] show that CTS can achieve significant reliability benefits by reducing the violations identified by RTCA, it is unclear how much economic benefits CTS can provide in RT SCED. Therefore, two EMS procedures, Procedure-A and Procedure-B, are proposed in this two-part paper to investigate the economic effects of implementing CTS in power system real-time operations.

This work was supported by 1) the Department of Energy Advanced Research Projects Agency - Energy, under the Green Electricity Network Integration program and under the Network Optimized Distributed Energy Systems program, and 2) the National Science Foundation award (1449080).

Xingpeng Li and Kory W. Hedman are with the School of Electrical, Computer, and Energy Engineering, Arizona State University, Tempe, AZ, 85287, USA (e-mail: Xingpeng.Li@asu.edu; kwh@myuw.net).

Procedure-A mimics the industrial practice which connects AC based RTCA and DC based RT SCED. Procedure-B enhances Procedure-A by integrating CTS into existing EMS with minimal change. Part-I of this paper includes a detailed literature review and presents the methodology, and Part-II introduces existing industrial practice and includes detailed results analysis.

The rest of this paper is organized as follows. Section II presents the industrial practice on RTCA, RT SCED, and CTS. Section III presents a comprehensive discussion of the results obtained with the traditional EMS Procedure-A and the enhanced EMS Procedure-B. The effects of integrating CTS into EMS on energy markets are described in Section IV. Section V presents the scalability studies of the proposed EMS procedures. Finally, Section VI concludes the paper.

## II. Industrial Practice

### A. Real-Time Contingency Analysis

To ensure secure operations of power systems, independent system operators (ISOs) must comply with reliability standards required by NERC. One important standard is $N$-1 reliability that requires power systems to withstand the loss of any single element. Thus, RTCA is conducted successively every few minutes at all ISOs. The actual implementation of contingency analysis could be different among different ISOs.

MISO's RTCA simulates more than 11,500 contingency scenarios every four minutes [5]. It is performed by solving contingency power flows independently. All the potential post-contingency flow violations and voltage violations are recorded, as well as the associated critical contingencies [6].

PJM runs a full AC contingency analysis to identify the contingencies that would cause system violations [7]. Approximately 6,000 contingencies are simulated every minute at PJM [7]. Though the PJM database has a list of all contingencies, not all of them are evaluated for each RTCA run [8].

A two-phase procedure is used in ERCOT to perform contingency analysis [9]. A heuristic screening is conducted in phase one to identify the most severe contingencies. In phase two, full AC analysis is performed on the contingencies identified in phase one and the contingencies specified in advance. ERCOT models about 4,000 contingencies in its database in 2012 and its RTCA engine executes every five minutes [10].

NYISO performs RTCA on pre-defined single and multiple contingencies. This would provide system operator with a list of potential transmission violations [11].

ISO-NE's RTCA runs every six minutes automatically or on demand [12]. It sorts violations by percent severity and provides operators with critical constraints that the associated post-contingency flow is over 90% of emergency limit [13].

CAISO's RTCA simulates about 2,200 pre-specified contingencies every five minutes [14]. It would report potential overloads and voltage violations under contingency, which alerts the operator to critical contingencies.

### B. Real-Time Security-Constrained Economic Dispatch

The classical AC OPF formulation was first developed in 1962 [15]. Although the problem has been formulated for over 50 years, a robust and reliable technique has not been developed to solve it in a timely manner due to the non-convexity and large-scale features. Thus, today's industry still uses the linearized DC power flow based SCED model.

RT SCED is a main application of the PJM real-time dispatch package. It does not change units' status and it dispatches online units in a single look-ahead period of 15 minutes. It runs about every five minutes or upon demand. For each run, three scenarios known as base scenario, high scenario, and low scenario are solved independently [16]. PJM's RT SCED co-optimizes energy, reserves, and regulation simultaneously [17]. It also provides a basis for the locational pricing calculator engine which runs exactly every five minutes to determine the locational marginal price (LMP) [17].

MISO's RT SCED dispatches energy and operating reserve to meet the forecasted demand and operating reserve requirements [18]. It executes continuously on a five-minute basis and the interval of its single look-ahead period is five minutes. It starts solving the problem five minutes before the target dispatch interval. Like PJM, MISO's SCED also uses a linear programming (LP) solver. Since commitment costs are sunk costs for SCED, the objective is to minimize the total operation cost that excludes start-up cost and no-load cost [18].

ISO-NE uses the unit dispatch system to perform SCED, which produces desired dispatch points (DDP) for the units in its territory. The DDP refreshes periodically on a five-minute basis [19] as the RT SCED runs every five minutes [20]. They must be approved by operators before sent to generators. The single-interval SCED jointly optimizes energy and reserves and typically looks 15 minutes ahead [19], [21]. It uses an incremental linear-optimization method to minimize the cost and produces DDP for dispatchable resources.

The real-time dispatch (RTD) application used by NYISO conducts SCED every five minutes with a look-ahead period of about an hour [22]. RTD is a multi-period dispatch process that simultaneously co-optimizes energy and reserves without involving commitment [23]. The RTD solution for the first five-minute interval would be implemented while the solutions for other intervals are for advisory purpose only. Real-time dispatch in corrective action mode executes on demand and overrides the regular RTD, and it may commit extra resources.

In ERCOT, RT SCED determines the least-cost dispatch solution to meet short-term load forecast. It executes every five minutes in ERCOT nodal market [24] and solves for a single interval of five minutes [25]-[26]. The SCED is a quadratic programming (QP) problem as the objective function is quadratic [25]-[26]. Currently, energy and ancillary services are not co-optimized in ERCOT real-time markets [27]-[28]. However, co-optimization is considered as a new initiative for improvement [27]. It produces LMPs and the price of system-wide reserves [29]. ERCOT runs SCED twice per cycle, which can reduce market power and ensure competition [30].

CAISO conducts SCED regularly on a five-minute basis and determines the least cost base points for participating units [31]. The SCED implemented at CAISO is a multi-period process that co-optimizes energy and ancillary services [32]. Only the solution of the first interval would be implemented. It also calculates LMP for market settlement.

Though the RT SCEDs implemented by various ISOs are very similar, they still have several different features. A comparison between various ISOs' SCEDs is presented in Table I. All ISOs run SCED every five minutes automatically or on demand. The SCEDs used by NYISO and CAISO look multi-interval ahead but only implement the first interval solutions; the other four ISOs implement single-interval SCED. All ISOs except for ERCOT co-optimize energy and reserves in their real-time markets. ERCOT's SCED is a QP problem while PJM, MISO, and ISO-NE execute SCED with LP solvers. To follow the most widely used features, the SCED implemented in this work is a LP based model that co-optimizes energy and reserves in a single interval of 15-minute.

Table I Comparison between various ISOs' RT SCED applications

| ISO | Execution cycle (minutes) | Type of periods | Only implement the solution of first period | Interval of the first period (minutes) | Look-ahead interval (minutes) | Co-optimize energy and reserve | Model |
|-----|------|------|------|------|------|------|------|
| PJM | 5 | single | NA | 15 | 15 | Yes | LP |
| MISO | 5 | single | NA | 5 | 5 | Yes | LP |
| ISO-NE | 5 | single | NA | 15 | 15 | Yes | LP |
| NYISO | 5 | multiple | Yes | 5 | ~60 | Yes | Unknown |
| CAISO | 5 | multiple | Yes | 5 | Unknown | Yes | Unknown |
| ERCOT | 5 | single | NA | 5 | 5 | No | QP |

NA denotes "not applicable", and Unknown means the associated information is not available publicly.

### C. Corrective Transmission Switching

Though studies on transmission switching in the literature started in 1980s [33]-[35], it has not been extensively used in industry today. The main concerns of implementing transmission switching include reduction of system security margin, instability issue due to discrete switching actions, and long computational time. However, with fast development of power engineering, optimization, and computer science technologies, those concerns and hurdles will be addressed eventually. Prior efforts in the literature have demonstrated that switching a line out of service does not necessarily adversely affect the system and can benefit the system in various aspects.

Transmission switching has gained a lot of attention recently. The hardware to implement it is circuit breaker that already exists in contemporary power systems. Switching a line out of service is fast enough for CTS to be considered as a realistic strategy in real-time operations. Thus, the requirement for actual implementation is to develop a decision support tool that can provide operators with beneficial CTS solutions.

As switching actions would degrade circuit breakers and reduce their lifespan, it would be practical if transmission switching is used as a corrective method or an emergency control scheme. Since the probability of a specific contingency is low, corrective transmission switching would rarely need to be implemented. Thus, the circuit breaker degradation due to CTS is negligible, which makes CTS a practical and promising control strategy for real-time operations.

In May 2009, due to the outages in the high voltage transmission system, significant congestion costs occurred for multiple days until CAISO was able to identify a switching action to relieve the congestion [36]. As documented in ISO-NE operating procedure [37], transmission switching is a viable option under both normal and emergency operating conditions and can be used to relieve transmission constraints. A list of switching solutions that serve as corrective actions in response to several specific contingencies is published by PJM [38]. In 2012, PJM took several high-voltage lines out of service as a corrective control in response to Superstorm Sandy to alleviate over-voltage problems [39].

Though there are several instances where transmission switching is implemented in practice to accomplish specific objectives, the decisions are predominantly made based on lookup table methods or ad hoc procedures that may require operators' personal judgment. Such empirical methods or offline analysis would limit the utilization of benefits provided by CTS. As a result, the implementation of CTS is limited. In addition, very little research work has examined how to practically utilize the reliability benefits provided by CTS in SCED. Therefore, accurate, fast, and systematic approaches for fully utilizing the network flexibility are essential for implementing CTS in industry.

## III. CASE STUDIES

In this section, the Cascadia system is used to verify the proposed Procedure-A and Procedure-B, as well as the proposed SCED models. This test case contains 179 buses, 245 branches, and 40 online units. The total online generation capacity is 9323 MW and the total load is 7324 MW. The Cascadia system is a synthetic test case based on the Washington state's geography; it is created and used by IncSys Academy for operator training purpose [40]. In this work, OpenPA [41] is used as the AC power flow solver and Gurobi is used as the optimization solver to solve SCED.

### A. Procedure-A: SCED with RTCA

To fully evaluate the proposed Procedure-A, base-case power flow is first performed; then, RTCA is conducted on all non-radial branch contingencies. No violation is observed in the base base. Fig. 1 shows the base-case condition of a key portion of the Cascadia system where contains two critical contingency-element and the beneficial switching actions.

The RTCA results on the Cascadia system are listed in Table II. Two out of 150 simulated contingencies cause overloads and thus, they are critical contingencies. These two criti-cal contingencies are parallel branches, branch 228 and branch 229. When one of them is out of service, the other branch would be overloaded by 241.6 MVA or 18.7% beyond the emergency rating. Fig. 2 shows branch 229 is overloaded under the outage of branch 228. As the two critical contingencies are equivalent and have the same consequences, only the contingency of branch 228 is analyzed for the rest of this paper.

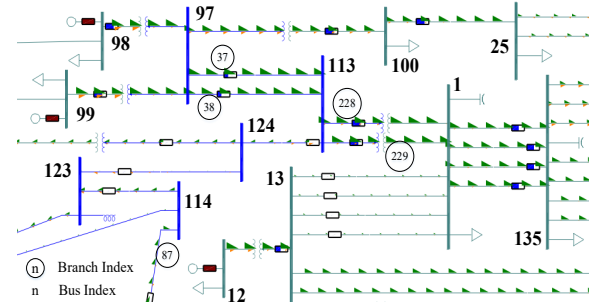

Fig. 1. System base-case condition of a portion of the Cascadia system

The proposed RT SCED models presented in Part-I of this paper are implemented and evaluated in this section. Thresholds $Pct$ and $PctC$ are set to 50% and 90% respectively unless they are explicitly described. The proposed SCED models share the majority of the constraints and the main difference between them is the network constraints.

Table III presents the results obtained with different SCED models. Cold-start PTDF power flow formulation based model SCED-M3 and cold-start $B$-$\theta$ power flow formulation based model SCED-M5 are essentially equivalent and they share the same lowest total cost among all SCED models. Moreover, the solution times for solving different SCED models are very similar. Hot-start PTDF based incremental model SCED-M1 results in the highest total cost, energy cost, and reserve cost.

Table II Results of RTCA on the Cascadia system

| Contingency branch | Monitor branch | Branch flow (MVA) | Emergency rating (MVA) | Violation (MVA) | Violation in percent |
|---|---|---|---|---|---|
| 228 | 229 | 1534.1 | 1292.5 | 241.6 | 18.7% |
| 229 | 228 | 1534.1 | 1292.5 | 241.6 | 18.7% |

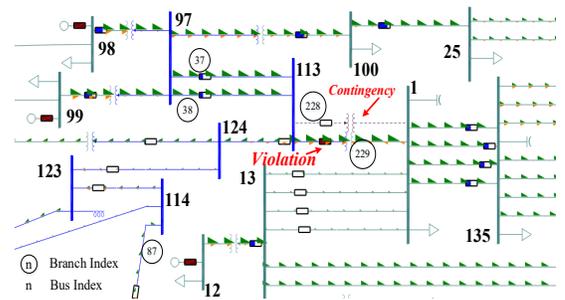

Fig. 2. System condition under the outage of branch 228.

Table III Results with different SCED models on the Cascadia system

| | Total cost ($/h) | Energy cost ($/h) | Reserve cost ($/h) | Solution time (s) |
|---|---|---|---|---|
| SCED-M1 | 50169.0 | 42943.3 | 7225.7 | 0.09 |
| SCED-M2 | 50011.8 | 42930.7 | 7081.1 | 0.11 |
| SCED-M3 | 49862.1 | 42903.7 | 6958.4 | 0.11 |
| SCED-M4 | 49903.8 | 42920.4 | 6983.3 | 0.14 |
| SCED-M5 | 49862.1 | 42903.7 | 6958.4 | 0.17 |

It is important to evaluate SCED solutions in an AC setting since DC power flow model used in SCED is an approximation and the accuracy is not guaranteed. Thus, the SCED solution, mainly the units' active power dispatch points, is fed back into base-case power flow simulation and contingency analysis simulation. The contingencies simulated in the post-SCED stage are the same with RTCA in the pre-SCED stage.

The results for SCED and post-SCED contingency analysis are shown in Table IV. Post-SCED RTCA causes two violations, 1) the overload on branch 229 under contingency of branch 228 and 2) the overload on branch 228 under contingency of branch 229. These two violations are equivalent as explained above and, thus, only the results for outage of branch 228 are listed in Table IV and other following tables.

As shown in Table IV, the results of RTCA with updated SCED solutions in the post-SCED situation illustrate the solution obtained with SCED-M1 outperforms the solutions determined by the other four models. With the dispatch points obtained by SCED-M1, the flow violation on branch 229 under outage of branch 228 is reduced from 241.6 MVA down to only 1.5 MVA, corresponding to an overload of 0.1% beyond the emergency limit, which is negligible; moreover, there are no other post-contingency violations or base-case violations. The dual variable of the network constraint on branch 229 under contingency of branch 228 with SCED-M1 is -10.5 $/MWh, which implies the total cost would be reduced by 10.5 $/h if the emergency limit increases by 1 MW.

Model SCED-M2 has the second-best performance. In this model, line outage distribution factor (LODF) is used to calculate the initial post-contingency branch flow. The solution of SCED-M2 causes a violation of 18.0 MVA for the same contingency of branch 228, which is 16.5 MVA more than that of SCED-M1. The extra 16.5 MVA overload comes from two sources: 1) DC model based LODF cannot accurately calculate the post-contingency branch flow, and 2) the branch emergency limit $LimitC_k$ calculated by (3) in Part-I is less precise than $LimitC_{kc}$ calculated by (2) in Part-I.

As expected, cold-start PTDF based SCED-M3 and cold-start $B$-$\theta$ based SCED-M5 share the least performance and result in 37.2 MVA violation on branch 229 under contingency 228. Model SCED-M4 has a better performance than SCED-M5 since the branch emergency limit used in SCED-M4 is more accurate; however, it still leads to a total violation of 29.8 MVA.

Though SCED-M1 results in the highest total cost as shown in Table III, it provides the best performance in the accurate AC setting. On the contrary, the dispatch solutions obtained by other SCED models would cause severe violations in an AC framework. In other words, SCED-M2 though SCED-M5 provide a cheaper solution at the cost of system security, which would violate security standards and put the system at a dangerous status. Therefore, the proposed model SCED-M1 is the preferred SCED model and the following evaluation and analysis for the rest part of this paper is based on SCED-M1.

Table IV Results of SCED and post-SCED contingency analysis with different SCED models on the Cascadia system

| | SCED | | | Post-SCED contingency analysis (branch 229 under contingency of branch 228) | | | |
|---|---|---|---|---|---|---|---|
| | Limit (MW) | Flow (MW) | Dual ($/MWh) | Rating (MVA) | Flow (MVA) | Violation (MVA) | Violation in percent |
| SCED-M1 | 1284.0 | 1284.0 | -10.5 | 1292.5 | 1294.0 | 1.5 | 0.1% |
| SCED-M2 | 1291.7 | 1291. | -9.8 | | 1310.6 | 18.0 | 1.4% |
| SCED-M3 | 1291.7 | 1291.7 | -2.4 | | 1329.7 | 37.2 | 2.9% |
| SCED-M4 | 1284.0 | 1284.0 | -9.8 | | 1322.4 | 29.8 | 2.3% |
| SCED-M5 | 1291.7 | 1291.7 | -2.4 | | 1329.7 | 37.2 | 2.9% |

Table V Results with different $Pct$ and $PctC$ on the Cascadia system

| $Pct$ | $PctC$ | SCED-M1 | | | Post-SCED contingency analysis (branch 229 under contingency of branch 228) | | | |
|---|---|---|---|---|---|---|---|---|
| | | Limit (MW) | Flow (MW) | Dual ($/MWh) | Rating (MVA) | Flow (MVA) | Violation (MVA) | Violation in percent |
| 1% | 1% | 1284.0 | 1284.0 | -10.5 | 1292.5 | 1294.0 | 1.5 | 0.1% |
| 50% | 50% | | 1284.0 | -10.5 | | 1294.0 | 1.5 | 0.1% |
| 80% | 80% | | 1284.0 | -10.5 | | 1294.0 | 1.5 | 0.1% |
| 100% | 100% | | 1284.0 | -10.5 | | 1294.0 | 1.5 | 0.1% |

As network constraints may largely affect SCED performance especially for large-scale real power systems, the sensitivity of thresholds $Pct$ and $PctC$ on the SCED performance is investigated in this work. The results for system performance are presented in Table V; in addition, the cost and computational time are presented in Table VI. It is worth noting that the system performance with different thresholds for selecting network constraints are the same, as well as the cost.

As shown in Table VI, the case with both $Pct$ and $PctC$ being set to 100% takes much less time than other cases while obtaining the same solutions. This is consistent with industrial practice. For real-time operations of a real power system, though the system condition changes, the initial dispatch point for SCED is not far away from the optimal solution, which is the key why only modeling a very small subset of critical network constraints can still maintain system security.

Table VI SCED results with different $Pct$ and $PctC$ on the Cascadia system

| $Pct$ | $PctC$ | SCED-M1 | | | |
|---|---|---|---|---|---|
| | | Total cost ($/h) | Energy cost ($/h) | Reserve cost ($/h) | Solution time (s) |
| 1% | 1% | 50169.0 | 42943.3 | 7225.7 | 5.27 |
| 50% | 50% | 50169.0 | 42943.3 | 7225.7 | 1.60 |
| 80% | 80% | 50169.0 | 42943.3 | 7225.7 | 0.12 |
| 100% | 100% | 50169.0 | 42943.3 | 7225.7 | 0.06 |

## B. Procedure-B: SCED with CTS-based RTCA

The proposed Procedure-B enhances Procedure-A by taking CTS into consideration. To focus on the potential benefits that can be provided by CTS, the duplicate results shared by both procedures are not presented again in this section.

In Procedure-B, CTS-based RTCA is implemented to augment the traditional RTCA. Table VII shows the CTS results for a contingency on branch 228. The top five best switching actions that provide Pareto improvement can reduce the post-contingency violation by 30.8%, 30.8%, 29.0%, 20.1%, and 19.9% respectively. Fig. 3 shows the system condition with branch 37 being switched off service for relieving overload. Though the overload on branch 229 still exists, it is reduced by 74.4 MVA with the top CTS solution.

Table VIII lists the emergency limits of branch 229 under contingency on branch 228 for different SCEDs. For a traditional SCED without CTS, the actual emergency limit is 1284.0 MW that is calculated by (2) in Part-I of this paper.

To take advantage of the violation reduction benefits provided by CTS, pseudo emergency limits are used in E-SCED to replace actual emergency limits. These pseudo emergency limits are higher than the actual emergency limits. E-SCED1, E-SCED2, E-SCED3, E-SCED4, and E-SCED5 use different pseudo emergency limits that are associated with the 1st, 2nd, 3rd, 4th, and 5th best switching actions respectively. The first best CTS solution can reduce the violation of the overloaded branch 229 by 74.4 MVA in the post-contingency situation with branch 228 being forced to be outage. This indicates that the emergency rating of branch 229 with CTS can be relaxed by 74.4 MVA under the contingency of branch 228. The associated pseudo emergency limit in MW for E-SCED can be calculated by (8) in Part-I of this two-part paper with the assumption that reactive power does not change in a short SCED period. Therefore, with the top CTS solution, the pseudo emergency limit of branch 229 under the contingency of branch 228 in E-SCED1 is 1358.8 MW, which is 74.8 MW higher than the actual emergency limit of 1284.0 MW. In addition, the emergency limit of branch 229 for E-SCED5 can increase by 48 MW even with the fifth best switching action.

To be consistent with the analysis in Section-III.A, $Pct$ and $PctC$ are set to be 50% and 90% respectively for E-SCED that considers CTS. Table IX presents the results of a traditional SCED without CTS, E-SCEDs with different CTS actions, and a relaxed SCED with no network constraint. It includes the results regarding cost, solution time, and dual variables of the network constraint for branch 229 under the contingency of branch 228 as well as its limit and post-SCED flow. It is worth noting that the size of the traditional SCED model is exactly the same with the enhanced SCED model. This shows that the enhanced SCED models with pseudo limits do not change the structure of traditional SCED models, which is also supported by their very similar solution time.

A binding network constraint may prevent cheap units from producing more power, which is the cause of congestion cost and unnecessary high total cost. The results of a SCED without consideration of network constraints are used as the benchmark to gauge the effects of CTS on SCED. By comparing the total costs of a traditional SCED and the SCED without any network constraints, the congestion cost of a traditional SCED without CTS can be calculated and it is 405.6 \$/h.

Table VII Results of RTCA with CTS on the Cascadia system

| Original violation on branch 229 (MVA) | CTS ranking | CTS branch | PI flag | Violation reduction (MVA) | Violation reduction in percent |
|---|---|---|---|---|---|
| 241.6 | 1st best | 37 | Yes | 74.4 | 30.8% |
| | 2nd best | 38 | Yes | 74.4 | 30.8% |
| | 3rd best | 85 | Yes | 70.1 | 29.0% |
| | 4th best | 86 | Yes | 48.6 | 20.1% |
| | 5th best | 87 | Yes | 48.0 | 19.9% |

Note that PI denotes Pareto improvement.

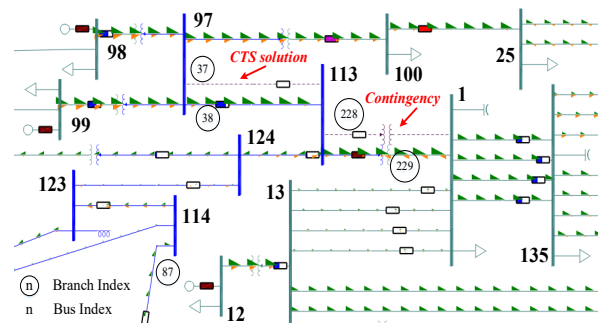

Fig. 3. System condition of a portion of the Cascadia system in the post-switching situation (CTS branch 37) under the outage of branch 228.

Table VIII Emergency limits of branch 229 under the contingency of branch 228 with and without CTS on the Cascadia system

| | CTS ranking | CTS branch | Actual emergency limit (in MW) w/o CTS | Pseudo emergency limit (in MW) w. CTS |
|---|---|---|---|---|
| SCED | NA | NA | 1284.0 | NA |
| E-SCED1 | 1st best | 37 | | 1358.8 |
| E-SCED2 | 2nd best | 38 | | 1358.8 |
| E-SCED3 | 3rd best | 85 | NA | 1354.5 |
| E-SCED4 | 4th best | 86 | | 1332.9 |
| E-SCED5 | 5th best | 87 | | 1332.3 |

NA denotes "not applicable".

Table IX Results of various SCEDs on the Cascadia system

| | | | Branch 229 under contingency 228 | | | Total cost (\$/h) | Congestion cost (\$/h) | Congestion cost reduction with CTS | LP Model | | | Solution time (s) |
|---|---|---|---|---|---|---|---|---|---|---|---|---|
| | | | Limit (MW) | Flow (MW) | Dual (\$/MWh) | | | | # of rows | # of columns | # of nonzeros | |
| E-SCED1 | with CTS | 1st best | 1358.8 | 1358.8 | -1.3 | 49797.9 | 34.5 | 91.5% | 6,526 | 5,243 | 35,687 | 0.11 |
| E-SCED2 | | 2nd best | 1358.8 | 1358.8 | -1.3 | 49797.9 | 34.5 | 91.5% | 6,526 | 5,243 | 35,687 | 0.11 |
| E-SCED3 | | 3rd best | 1354.5 | 1354.5 | -1.4 | 49803.6 | 40.2 | 90.1% | 6,526 | 5,243 | 35,687 | 0.12 |
| E-SCED4 | | 4th best | 1332.9 | 1332.9 | -1.8 | 49834.6 | 71.2 | 82.5% | 6,526 | 5,243 | 35,687 | 0.12 |
| E-SCED5 | | 5th best | 1332.3 | 1332.3 | -1.8 | 49835.8 | 71.6 | 82.4% | 6,526 | 5,243 | 35,687 | 0.12 |
| SCED | w/o. CTS | | 1284.0 | 1284.0 | -10.5 | 50169.0 | 405.6 | NA | 6,526 | 5,243 | 35,687 | 0.09 |
| Relaxed SCED | With no network constraint | | NA | | | 49763.4 | 0.0 | NA | 1,126 | 1,001 | 4,346 | 0.04 |

NA denotes "not applicable".

Table IX illustrates that congestion cost drops with a higher limit of the bottleneck branch. The congestion cost of SCED is reduced by 91.5% with the top identified switching actions or 82.4% with the fifth best identified switching actions. This demonstrates that using pseudo limit in E-SCED can achieve significant congestion cost reduction and thus improve the social welfare. Fig. 4 presents the congestion costs of a traditional SCED and multiple E-SCEDs. With the top five identified CTS solutions being considered in E-SCEDs, the congestion cost is reduced from 405.6 $/h to a much smaller value ranging from 34.5 $/h to 71.6 $/h.

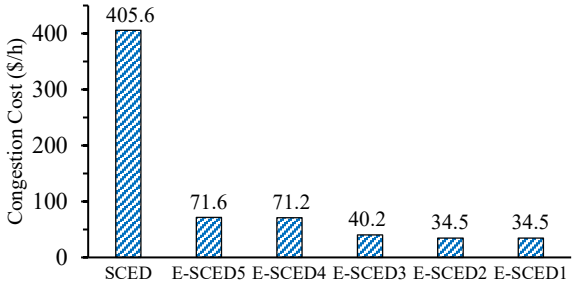

Fig. 4. Congestion costs of the traditional SCED and various E-SCEDs on the Cascadia system.

It is worth noting that the top switching action can reduce the congestion cost by 91.5% while it can only reduce the associated post-contingency violation by 30.8%. By implementing CTS, the congestion cost reduction in percent is typically higher than post-contingency violation reduction in percent. The reason is that the marginal cost reduction due to 1 MW increase of branch limit may drop as the associated branch limit becomes higher. Relieving binding network constraint by increasing branch limit would allow the cheapest available unit to ramp up and force the most expensive unit to reduce its output. If that branch is still binding after the cheapest available unit reaches its maximum and the most expensive unit reaches its minimum, further increasing its limit may allow the second cheapest available unit to ramp up and reduce the output of the second most expensive unit. This would still reduce the total cost, but the cost reduction would drop for each MW relieved in the branch limit as the limit increases.

The above conclusion can also be made from the dual variables of network constraints. With the actual emergency limit 1284.0 MW of branch 229 under contingency 228, the associated dual variable is -10.5 $/h, which implies that the total cost would drop by 10.5 $/h if its emergency limit increases by 1 MW from 1284.0 MW. However, when the pseudo emergency limit 1332.3 MW is used, the associated dual variable becomes -1.8 $/h, which implies the total cost would only drop by 1.8 $/h if the emergency limit increases by 1 MW from 1332.3 MW. Thus, as the branch limit increases, the marginal cost reduction drops, which indicates CTS can achieve higher cost reduction in percent than violation reduction in percent.

Though it has been demonstrated that congestion cost can be significantly reduced with CTS by using pseudo emergency limits in E-SCED, the violation reduction performance of CTS in the post-E-SCED stage should be examined as the system condition changes.

In the pre-E-SCED stage, post-contingency violation reductions are different with different CTS solutions and, thus, the associated pseudo emergency limits are also different in the E-SCED models, which may result in different dispatch solutions. The results of traditional RTCA in the post-E-SCED stage, with different CTS or pseudo limits considered in E-SCED, are shown in Table X. The results in this table are for branch 229 under contingency 228. As expected, the post-contingency violation increases with higher ranked beneficial switching action considered in E-SCED.

Table X Results of traditional RTCA in the post-E-SCED situations

| CTS for E-SCED | Actual emergency rating (MVA) | Flow (MVA) | Violation (MVA) | Violation (%) |
|---|---|---|---|---|
| 1st best | | | | |
| 2nd best | | 1376.5 | 84.0 | 6.50% |
| 3rd best | 1292.5 | 1371.6 | 79.0 | 6.11% |
| 4th best | | 1346.6 | 54.1 | 4.18% |
| 5th best | | 1345.8 | 53.3 | 4.12% |

Table XI Results of RTCA with CTS in the post-E-SCED stage with the SCED solution corresponding to the 1st best and 2nd best CTS solutions identified in the pre-E-SCED stage on the Cascadia system

| CTS for RTCA | Flow (MVA) | Flow change caused by CTS (MVA) | Violation (MVA) | Violation in percent (%) | Violation reduction (MVA) | Violation reduction in percent |
|---|---|---|---|---|---|---|
| 1st best | 1304.9 | -71.5 | 12.4 | 0.96% | 71.5 | 85.2% |
| 2nd best | 1304.9 | -71.5 | 12.4 | 0.96% | 71.5 | 85.2% |
| 3rd best | 1235.6 | -140.9 | -56.9 | 0.0% | 84.0 | 100% |
| 4th best | 1285.8 | -90.7 | -6.7 | 0.0% | 84.0 | 100% |
| 5th best | 1280.2 | -96.3 | -12.3 | 0.0% | 84.0 | 100% |

Table XI shows the results of CTS-based RTCA in the post-E-SCED stage. The E-SCED solution used for Table XI corresponds to the 1st best and 2nd best CTS solution identified in the pre-E-SCED stage. Note that the 1st best CTS branch 37 and the 2nd best CTS branch 38 are equivalent since they are in parallel and have the same parameters. The results in this table are for branch 229 under contingency 228. Before implementing CTS in the post-E-SCED stage, contingency 228 causes a violation of 84 MVA on branch 229. However, that violation can be relieved with the five beneficial CTS solutions identified in the pre-E-SCED stage. The top two switching actions can reduce the violation by about 85% while the other three switching actions can fully eliminate the violation.

Though the amounts of violation reduction with CTS in the post-E-SCED scenario are slightly different with pre-E-SCED scenario, all CTS actions identified in the pre-E-SCED scenario can reduce flow on the same overloaded branch in the post-E-SCED stage. This demonstrates that CTS can provide benefits even when the system condition changes.

Table XII, XIII, and Table XIV show the results of CTS-based RTCA in the post-E-SCED stage, where E-SCEDs correspond to the 3rd, 4th, and 5th best CTS solutions respectively. If the pseudo emergency limit associated with the 3rd best CTS

action is used in E-SCED, RTCA simulated in the post-E-SCED stage results in an overload of 79 MVA on branch 229 under contingency 228; however, three of the five CTS actions can fully eliminate the post-contingency violation while the other two CTS solutions can reduce the overload by more than 90%. For E-SCED with the 4th or 5th best CTS action, all five CTS actions can fully eliminate the overload.

With lower limits being used for the network constraints in E-SCED, branches would have more security margins in the post-E-SCED stage. As the congestion cost reduction does not vary much with different CTS solutions considered in E-SCED, using the pseudo limit associated with the 3rd best switching actions can provide both substantial economic benefits and significant post-contingency violation reductions.

Table XII Results of RTCA with CTS in the post-E-SCED stage with the SCED solution corresponding to the 3rd best CTS solution identified in the pre-E-SCED stage on the Cascadia system

| CTS for RTCA | Flow (MVA) | Violation (%) | Violation (MVA) | Violation reduction (MVA) | Violation reduction (%) |
|---|---|---|---|---|---|
| 1st best | 1300.0 | 0.6% | 7.5 | 71.5 | 90.5% |
| 2nd best | 1300.0 | 0.6% | 7.5 | 71.5 | 90.5% |
| 3rd best | 1228.5 | 0.0% | -64.1 | 79.0 | 100% |
| 4th best | 1282.2 | 0.0% | -10.3 | 79.0 | 100% |
| 5th best | 1277.0 | 0.0% | -15.5 | 79.0 | 100% |

Table XIII Results of RTCA with CTS in the post-E-SCED stage with the SCED solution corresponding to the 4th best CTS solution identified in the pre-E-SCED stage on the Cascadia system

| CTS | Flow (MVA) | Violation (MVA) | Violation (%) | Violation reduction (MVA) | Violation reduction (%) |
|---|---|---|---|---|---|
| 1st best | 1275.3 | -17.2 | 0.0% | 54.1 | 100% |
| 2nd best | 1275.3 | -17.2 | 0.0% | 54.1 | 100% |
| 3rd best | 1195.2 | -97.3 | 0.0% | 54.1 | 100% |
| 4th best | 1263.8 | -28.7 | 0.0% | 54.1 | 100% |
| 5th best | 1260.8 | -31.7 | 0.0% | 54.1 | 100% |

Table XIV Results of RTCA with CTS in the post-E-SCED stage with the SCED solution corresponding to the 5th best CTS solution identified in the pre-E-SCED stage on the Cascadia system

| CTS | Flow (MVA) | Violation (MVA) | Violation (%) | Violation reduction (MVA) | Violation reduction (%) |
|---|---|---|---|---|---|
| 1st best | 1274.5 | -18 | 0.0% | 53.3 | 100% |
| 2nd best | 1274.5 | -18 | 0.0% | 53.3 | 100% |
| 3rd best | 1194.4 | -98.1 | 0.0% | 53.3 | 100% |
| 4th best | 1263.1 | -29.4 | 0.0% | 53.3 | 100% |
| 5th best | 1260.2 | -32.3 | 0.0% | 53.3 | 100% |

## IV. MARKET IMPLICATION

Market results of the traditional SCED and the proposed E-SCEDs are presented in Table XV. When the flexibility in transmission network is taken into account, the load payment drops significantly, as well as the generator revenue, generator rent, and congestion revenue. It is observed that with CTS, the amount of load payment reduction is much more than the amount of generator rent reduction. Fig. 5 and Fig. 6 show the load payment and congestion revenue respectively, for the traditional SCED and the proposed E-SCED. Apparently, with higher pseudo limit used in E-SCED, the system-wide load payment and congestion revenue substantially decrease, which implies that the proposed E-SCED can improve the market efficiency in comparison with a traditional SCED.

The nodal LMP including energy component and congestion component is shown in Table XVI. The energy LMP of each bus is the same across the entire system and is also equal to the LMP at the slack bus. The average LMPs and the average congestion LMPs for various E-SCEDs are very close, since even the fifth best CTS solution can relieve the congestion by 82%. In comparison with the traditional SCED, the average LMP is reduced by about 5%; and the average congestion LMP is reduced by 82% to 88%, which is consistent with the degree of congestion cost reduction in percent as shown in Table IX. This shows the application of pseudo limit in optimal generation dispatching can reduce LMP by a significant level, which explains the results shown in Table XV.

Table XV Market results with a traditional SCED and various E-SCED

| | Load payment ($/h) | Generator revenue ($/h) | Generator cost ($/h) | Generator rent ($/h) | Congestion revenue ($/h) |
|---|---|---|---|---|---|
| E-SCED1 | 58158.5 | 57291.0 | 42839.5 | 14451.5 | 867.5 |
| E-SCED2 | 58158.5 | 57291.0 | 42839.5 | 14451.5 | 867.5 |
| E-SCED3 | 58112.1 | 57056.2 | 42845.2 | 14211.0 | 1055.9 |
| E-SCED4 | 58977.2 | 57363.6 | 42876.2 | 14487.4 | 1613.6 |
| E-SCED5 | 58977.2 | 57364.7 | 42877.3 | 14487.4 | 1612.5 |
| SCED | 74865.8 | 62553.6 | 42943.3 | 19610.3 | 12312.2 |

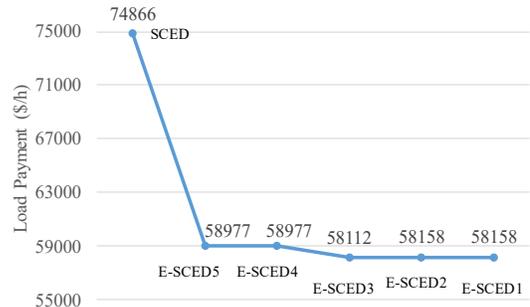

Fig. 5. Load payment for SCED and E-SCEDs on the Cascadia system.

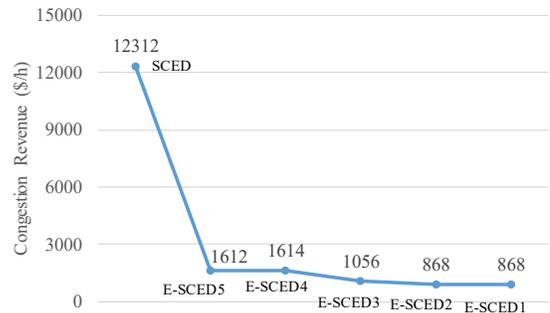

Fig. 6. Congestion revenue for SCED and E-SCEDs on the Cascadia system.

Table XVI Average LMP with SCED and various E-SCEDs on Cascadia

| | Average LMP ($/MWh) | Average congestion LMP ($/MWh) | Energy LMP ($/MWh) |
|---|---|---|---|
| E-SCED1 | 7.685 | 0.063 | 7.622 |
| E-SCED2 | 7.685 | 0.063 | 7.622 |
| E-SCED3 | 7.649 | 0.071 | 7.578 |
| E-SCED4 | 7.672 | 0.094 | 7.578 |
| E-SCED5 | 7.672 | 0.095 | 7.578 |
| SCED | 8.037 | 0.542 | 7.494 |

## V. SCALABILITY STUDIES

To better validate the proposed EMS procedures, scalability studies of both Procedure-A and Procedure-B are performed on the large-scale Polish system in this section. This system has 2,382 buses, 2,895 branches, and 327 generators. Its total amount of load is 21.5 GW.

**Pre-SCED Phase**

In the pre-SCED phase, base-case AC power flow simulation reports a total violation of 15.1 MVA on two lines; and then AC RTCA conducts and identifies 33 critical contingencies causing 46 potential violations with a total amount of 434.7 MVA, which corresponds to 13.2 MVA violation per contingency or 9.5 MVA per violation. For Procedure-B, CTS is conducted to determine beneficial switching lines that can reduce the post-contingency overloads and the results with Pareto improvement are shown in Table XVII. The average violation reduction in percent [1] for top five best switching solutions, denoted by $\tau_{CTS}$, ranges from 74.9% to 22.7%. This indicates that corrective transmission switching can effectively relieve flow violations on the Polish system.

Table XVII Results of RTCA with CTS on the Polish system

| CTS Ranking | 1st best | 2nd best | 3rd best | 4th best | 5th best |
|---|---|---|---|---|---|
| $\tau_{CTS}$ | 74.9% | 62.2% | 45.1% | 28.2% | 22.7% |

**SCED Phase**

In this SCED phase, both procedures are implemented to obtain the least cost solutions while addressing the network violations. The traditional SCED of Procedure-A uses the network constraints reported from RTCA while the E-SCED of Procedure-B uses the relaxed network constraints with pseudo limits determined from the results of RTCA with 3rd best CTS solutions. To be consistent, the thresholds $Pct$ and $PctC$ for both procedures are set to 50% and 90% respectively for this scalability studies. In addition, a relaxed SCED without any network constraints is also solved to provide a benchmark to gauge the proposed procedures.

The detailed results of the three SCED models are presented in Table XVIII. It is observed that the congestion cost for the traditional SCED of Procedure-A is 14,328 $/h, which is reduced to 9,257 $/h for E-SCED of Procedure-B. In other words, the use of branch pseudo emergency limits in the application of SCED can reduce the system congestion cost by 5,071 $/h or a level of 35.4% for the large-scale Polish system.

Table XVIII Results of various SCED models on the Polish system

|  |  | Relaxed SCED with no network | Traditional SCED | E-SCED with the 3rd best CTS |
|---|---|---|---|---|
| Total cost ($/h) | | 262,006 | 276,337 | 271,266 |
| Congestion cost ($/h) | | NA | 14,328 | 9257 |
| Congestion cost reduction with CTS ($/h) | | NA | NA | 5,071 (35.4%) |
| Average LMP ($/MWh) | | 16.7 | 28.5 | 19.9 |
| LP model size | # of rows | 12,478 | 22,157 | 22,157 |
| | # of columns | 11,492 | 13,862 | 13,862 |
| | # of nonzeros | 127,676 | 19,980,218 | 19,980,218 |
| # of iterations taken | | 233 | 2,612 | 2,614 |
| Solution time (s) | | 0.12 | 9.69 | 9.62 |

NA denotes "not applicable".

Moreover, as the LP model statistics in Table XVIII shows, the size of the traditional SCED model is exactly the same with the size of the E-SCED model using branch pseudo emergency limits; the solution time and number of LP iterations used for both problems to converge are very similar. This demonstrates that the proposed E-SCED model with pseudo emergency limits does not increase the problem complexity and can solve in a very similar time. The number of nonzeros for the traditional SCED and E-SCED models on the large-scale Polish system is 19,980,218, which is 500 times larger than the Cascadia system that has 35,687 nonzeros. The system size increases by about 13 times in terms of number of buses while the nonzeros of the SCED LP problem increases by more than 500 times; however, the solution time only increases by a factor of less than 100. This indicates that the proposed EMS procedures can scale very well for large-scale power systems.

When this is no network congestion, the nodal LMP is the same, 16.7 $/MWh, throughout the entire network; and the load payment is 360,438 $/h. With the traditional SCED model using the RTAC information directly, the system average LMP becomes 28.5 $/MWh and the load payment increases to 695,484 $/h; the generator revenue is 520,798 $/h and thus, the network congestion revenue is 174,686 $/h. These market numbers can go down significantly with the use of pseudo limit in E-SCED. With E-SCED, the system average LMP, load payment, generator revenue, and congestion revenue drop to 19.9 $/MWh, 457,047 $/h, 391,682 $/h, and 65,365 $/h. Apparently, the E-SCED model of the proposed Procedure-B can significantly increase the social welfare including a reduction of 62.6% in the system congestion revenue.

Table XIX shows the statistics of the network constraints for the traditional SCED model and the E-SCED models. The number of base-case branch constraints and contingency-case branch constraints are the same. They also share the same 2 critical base-case branch constraints and 46 critical contingency-case branch constraints that are related to the violations identified by base-case power flow simulation and RTCA respectively. With E-SCED, the emergency limits of 26 critical contingency-case branch constraints can be relaxed and replaced with the proposed pseudo emergency limits; the amount of increased branch emergency limit is 2.8 MW on average. Obviously, these 26 relaxed critical contingency-case branch constraints account for the significant economic benefits achieved with E-SCED.

Table XIX Statistics of network constraints for SCED and E-SCED

|  | Traditional SCED | E-SCED with the 3rd best CTS |
|---|---|---|
| # of all branch constraints | 2539 | 2539 |
| # of base-case branch constraints | 169 | 1639 |
| # of contingency-case branch constraints | 2,370 | 2,370 |
| # of critical branch constraints related to violations | 48 | 48 |
| # of critical contingency-case branch constraints with pseudo limits | NA | 26 |

NA denotes "not applicable".

**Post-SCED Phase**

The achievement of violation elimination, congestion relief, and cost reduction cannot compromise the system reliability. Thus, to validate the traditional SCED of the proposed Procedure-A and the E-SCED of the proposed Procedure-B, RTCA and CTS-based RTCA are performed on the large-scale Polish system with updated generations in the post-SCED phase.

With the generator dispatch points determined by the traditional SCED of Procedure-A, the 15.1 MVA actual base-case violations in the pre-SCED phase are fully eliminated and the 434.7 MVA potential post-contingency violations are reduced to only 0.9 MVA that is only 0.2% of the post-contingency violations identified in the pre-SCED phase. This demonstrates the proposed Procedure-A can effectively eliminate network violations. Similarly, with the generator dispatch points determined by E-SCED of Procedure-B, the 15.1 MVA actual base-case violations are also fully eliminated and the 434.7 MVA potential post-contingency violations are reduced to 6.3 MVA that is 1.4% of the original post-contingency violations. Five potential post-contingency violations exist in the post-E-SCED phase and the largest overload is 4.3 MVA, which is not trivial; this is expected since some branch limits are relaxed in the E-SCED model. However, all these five violations can be fully eliminated with the same CTS solutions (multiple switching actions available) identified in the pre-E-SCED phase. Therefore, the RTCA with CTS studies strongly support the proposed EMS Procedure-B using branch pseudo emergency limit in the E-SCED model.

## VI. Conclusions

A novel EMS procedure, Procedure-B, is proposed in this paper for control centers to take advantage of the flexibility in transmission networks by means of corrective transmission switching in the real-time operations of electric power systems. With the proposed EMS procedure, the reliability benefits provided by CTS in RTCA can be forwarded to SCED for economic benefits. With the proposed pseudo limit, Procedure-B requires minimal change in existing EMSs, which will minimize the obstacles of implementing the proposed Procedure-B by simply adding one separate CTS module after RTCA and before SCED.

Numerical simulations conducted in Part-II of this paper demonstrate the effectiveness of considering the network as a flexible asset in power system real-time operations. The results show the proposed Procedure-B can take advantage of the flexibility in the transmission networks. With Procedure-B, corrective transmission switching can reduce substantial post-contingency violations identified in RTCA and achieve significant economic benefits in E-SCED. The system congestion cost can be largely reduced with E-SCED by using the proposed concept of pseudo limit; in addition, the proposed E-SCED does not increase the problem size and can solve in a very similar timeframe of the traditional SCED. The analysis based on comparing the real-time market results of the traditional SCED and the proposed E-SCED shows that E-SCED with pseudo limit can significantly improve the social welfare by reducing the system overall LMP, total load payment, and total operation cost. The scalability studies demonstrate the proposed EMS procedures can scale well for large-scale power systems.

**Xingpeng Li** (S'13−M'18) received the B.S. degree in electrical engineering from Shandong University, Jinan, China, in 2010, and the M.S. degree in electrical engineering from Zhejiang University, Hangzhou, China, in 2013. He received the M.S. degree in industrial engineering and the Ph.D. degree in electrical engineering from Arizona State University, Tempe, AZ, USA, in 2016 and 2017 respectively.

Currently, he is an Assistant Professor in the Department of Electrical and Computer Engineering at the University of Houston. He previously worked for ISO New England, Holyoke, MA, USA, and PJM Interconnection, Audubon, PA, USA. Before joining the University of Houston, he was a Senior Application Engineer for ABB, San Jose, CA, USA. His research interests include power system operations and optimization, energy management system, energy markets, microgrids, and grid integration of renewable energy sources.

**Kory W. Hedman** (S'05−M'10−SM'16) received a B.S. degree in electrical engineering and a B.S. degree in economics from the University of Washington, Seattle, WA, USA, in 2004, an M.S. degree in economics and an M.S. degree in electrical engineering from Iowa State University, Ames, IA, USA, in 2006 and 2007, respectively, and the M.S. and Ph.D. degrees in operations research from the University of California—Berkeley, Berkeley, CA, USA, in 2008 and 2010, respectively.

He is currently an Associate Professor with the School of Electrical, Computer, and Energy Engineering, Arizona State University, Tempe, AZ, USA. His research interests include power systems operations and planning, market design, power system economics, renewable energy, and operations research. He was the recipient of the Presidential Early Career Award for Scientists and Engineers from the U.S. President Barack H. Obama in 2017.